\documentclass[11pt]{article}
\usepackage{amssymb,amsfonts,amsmath,latexsym, epsfig,mathrsfs,colordvi}
\usepackage{xcolor,graphicx,float}
\usepackage{tikz}
\usepackage{graphicx}
\newtheorem{theo}{Theorem}

\newtheorem{exam}[theo]{Example}
\newtheorem{lem} [theo]{Lemma}
\newtheorem{coro}[theo]{Corollary}

\makeatletter \@addtoreset{equation}{section}
\@addtoreset{theo}{}\makeatother

\def\qed{\hfill \rule{4pt}{7pt}}
\def\pf{\noindent {\it Proof.} }

\def\Peaks{\mathrm{Peaks}}
\def\Humps{\mathrm{Humps}}

\begin{document}

\title{Enumerations of humps and peaks in $(k,a)$-paths and $(n,m)$-Dyck paths
via bijective proofs}

\author{Rosena R. X. Du\footnote{Corresponding Author. Email: rxdu@math.ecnu.edu.cn.}, Yingying Nie and Xuezhi Sun\\ \\ Department of Mathematics, Shanghai Key Laboratory of PMMP \\East China Normal University,
500 Dongchuan Road \\Shanghai, 200241, P. R. China.}

\date{March 26, 2015}
\maketitle

\vskip 0.7cm \noindent{\bf Abstract.}
Recently Mansour and Shattuck studied $(k,a)$-paths and gave formulas that related the total number of humps in all $(k,a)$-paths to the number of super $(k,a)$-paths. These results generalized earlier results of Regev on Dyck paths and Motzkin paths. Their proofs are based on generating functions and they asked for bijective proofs for their results. In this paper we first give bijective proofs of Mansour and Shattuck's results, then we extend our study to $(n,m)$-Dyck paths. We give a bijection that relates the total number of peaks in all $(n,m)$-Dyck paths to certain free $(n,m)$-paths when $n$ and $m$ are coprime. From this bijection we get the number of $(n,m)$-Dyck paths with exactly $j$ peaks, which is a generalization of the well-known result that the number Dyck paths of order $n$ with exactly $j$ peaks is the Narayana number $\frac{1}{k}{n-1\choose k-1}{n\choose k-1}$.



\section{Introduction}

In this paper we study two kinds of lattice paths: $(k,a)$-paths and $(n,m)$-Dyck paths.

A $(k,a)$-path of order $n\in \mathbb{N}$ where $k,a \in \mathbb{Z}^+$ is a lattice path in $\mathbb{Z}\times\mathbb{Z}$ from $(0,0)$ to $(n,0)$ which uses up steps $(1,k)$, down steps $(1,-1)$ and horizontal steps $(a,0)$ and never goes below the $x$-axis. We use $\mathcal{P}_n(k,a)$ to denote the set of all $(k,a)$-paths of order $n$. Note that $\mathcal{P}_n(1,\infty)$, $\mathcal{P}_n(1,1)$ and $P_n(1,2)$ are the set of Dyck paths, Motzkin paths, and Schr\"{o}der paths, respectively. And $\mathcal{P}_n(k, \infty)$ denotes the set of $(k,a)$-paths of order $n$ without horizontal steps. Such paths are also called $k$-ary paths and are studied in \cite{LiMansour}.

An $(n,m)$-\emph{Dyck path} where $n,m \in \mathbb{Z}^+$ is a lattice path in $\mathbb{Z}\times\mathbb{Z}$, from $(0,0)$ to $(n,m)$, which uses up steps $(0,1)$ and down steps $(1,0)$ and never goes below the diagonal line $y=\frac{m}{n}x$.
Note that when $m=n$, if we rotate an $(n,n)$-Dyck path 45 degrees clockwise, we get an ordinary Dyck path of order $n$. Here we prefer to use this slightly different form of definition for later convenience. (This is why we call the step $(1,0)$ a ``down" step instead of a ``right" step).

$(n,m)$-Dyck paths have been studied by many authors \cite{DArmstrong,Bizley,Y.FUKUKAWA}. It is known that when $\gcd(n,m)=1$, i.e., when $n$ and $m$ are coprime, the number of $(n,m)$-Dyck paths is:
\begin{equation}
D(n,m)=\frac{1}{n+m}{n+m\choose n}.
\end{equation}

In this paper we focus on counting these paths with a given number of peaks or humps.
A {\em peak} in a path is an up step followed immediately by a down step. A {\em hump} in a $(k,a)$-path is an up step followed by zero or more horizontal steps followed by a down step. We denote by $\#\Peaks(P)$ ($\#\Humps(P)$) the number of peaks (humps) in a path $P$. If a $(k,a)$-path is allowed to go below the $x$-axis, then we call it a {\em super} path, or a {\em free} path. Let $\mathcal{SP}_n(k,a)$ denote the set of all super $(k,a)$-paths of order $n$.
In \cite{HumpsRegev} Regev noticed the following curious relation between the number of peaks in Dyck paths and the number of free paths:

\begin{equation}\label{peak1}
2\sum_{P\in \mathcal{P}_n(1,\infty)}\#\Peaks(P)=|\mathcal{SP}_{n}(1,\infty)|.
\end{equation}

He also proved that the following equation holds for the number of humps in Motzkin paths and super Motzkin paths.
\begin{equation}\label{hump1}
2\sum_{P\in \mathcal{P}_n(1,1)}\#\Humps(P)=|\mathcal{SP}_{n}(1,1)|-1.
\end{equation}

Regev's proofs of the above two equations involve a recurrence relation and the WZ method \cite{A=B,Zeilberger}. And he asked for bijective proofs for these equations in \cite{HumpsRegev}. In \cite{DingDu} Ding and Du gave bijective proofs of these two equations, and also proved that similar relations hold for Schr\"{o}der paths.

Recently, using generating function methods, Mansour and Shattuck \cite{MansourShattuck} generalized the above results to ordinary $(k,a)$-paths and proved the following equations:
\begin{equation}\label{hump}
(k+1)\sum_{P\in \mathcal{P}_n(k,a)}\#\Humps(P)=|\mathcal{SP}_{n}(k,a)|-\delta_{a|n},
\end{equation}

\begin{equation}\label{peak}
(k+1)\sum_{P\in \mathcal{P}_n(k,a)}\#\Peaks(P)=|\mathcal{SP}_{n}(k,a)|-|\mathcal{SP}_{n-a}(k,a)|,
\end{equation}
where $\delta_{a|n}=1$ if $a$ divides $n$ or $0$ otherwise. Putting $k=1$ and $a=1$, $a=\infty$ in \eqref{hump} and \eqref{peak} we get \eqref{hump1} and \eqref{peak1}.

%
%

The outline of this paper is as follows. In Section 2 we define a bijection and prove Equations \eqref{hump} and \eqref{peak} for $(k,a)$-paths. In Section 3 we first study the properties of $(n,m)$-Dyck paths when $\gcd(n,m)=1$. Then we give a bijection and count the number of $(n,m)$-Dyck paths with a given number of peaks. In Section 4 we discuss the special case of $(n,kn)$-Dyck paths, which are in one-to-one correspondence with $k$-ary paths, and give a formula for the number of $k$-ary paths with a given number of peaks.

\section{Bijective Proofs for Mansour and Shattuck's Results}



Observe that equation \eqref{hump} can also be written in the following form:
\begin{equation}\label{1overk+1}
\sum_{P\in\mathcal{P}_{n}(k,a)}\#\Humps(P)=\frac{1}{k+1}\left(|\mathcal {SP}_{n}(k,a)|-\delta_{a|n}\right).
\end{equation}
Instead of proving \eqref{hump} directly, we will give a bijection which proves \eqref{1overk+1}.

Let $\mathcal{SP}_{n}^{0}(k,a)$ denote the set of super $(k,a)$-paths in $\mathcal{SP}_{n}(k,a)$ that contains at least one up step, and let $\mathcal{SP}_{n}^{U}(k,a)$ denote the set of super  $(k,a)$-paths in $\mathcal{SP}^{0}_{n}(k,a)$ whose first non-horizontal step is an up step. Note that there is one super $(k,a)$-path in $\mathcal {SP}_{n}(k,a)$ that consists of only horizontal steps if and only if $n$ is divisible by $a$. Hence $|\mathcal{SP}_n^0(k,a)|=|\mathcal{SP}_n(k,a)|-\delta_{a\mid n}$. The following lemma explains what the right hand side of \eqref{1overk+1} counts.

\begin{figure}[h]
\begin{center}
\begin{tikzpicture}[mark=*,mark size=0.2ex,scale=0.6]
\draw plot coordinates %
  {(0,0) (10pt,0pt)(20pt,0pt)(30pt,0pt)(40pt,0pt)};
  \node at(20pt,-10pt){$H^{l}$};

\draw[very thick,color=red] plot coordinates %
  {(40pt,0pt)(50pt,50pt)};
  \draw[xshift=50pt,yshift=50pt] (0pt,0pt) .. controls (5pt,40pt)%
and (20pt,20pt) .. (40pt,0);
\draw[xshift=50pt,yshift=50pt,very thick,color=blue] plot coordinates %
  {(40pt,0pt)(50pt,-10pt)};
 \draw[dashed,mark=,xshift=50pt,yshift=50pt] plot coordinates %
      {(0pt,0pt)(40pt,0pt)};

      \draw[xshift=128pt,yshift=30pt,very thick, color=blue] plot coordinates %
  {(-10pt,10pt)(0pt,0pt)};
       \draw[xshift=128pt,yshift=30pt] (0pt,0pt) .. controls (5pt,40pt)%
and (20pt,20pt) .. (40pt,0);
  \draw[dashed,mark=,xshift=128pt,yshift=30pt] plot coordinates %
      {(0pt,0pt)(40pt,0pt)};

 \draw[xshift=128pt,yshift=30pt,very thick, color=blue] plot coordinates %
  {(40pt,0pt)(50pt,-10pt)};
\draw[xshift=178pt,yshift=20pt,very thick, color=blue] plot coordinates %
  {(40pt,0pt)(50pt,-10pt)};
  \draw[xshift=178pt,yshift=20pt] (0pt,0pt) .. controls (5pt,40pt)%
and (20pt,20pt) .. (40pt,0);
  \draw[dashed,mark=,xshift=178pt,yshift=20pt] plot coordinates %
      {(0pt,0pt)(40pt,0pt)};
\draw[dashed]plot coordinates %
      {(0pt,0pt)(40pt,0pt)};

      \node at(109pt,40pt){$\ldots$};

 \draw[xshift=205pt,yshift=10pt,very thick, color=blue] plot coordinates %
  {(40pt,0pt)(50pt,-10pt)};
    \draw[xshift=255pt,yshift=0pt] (0pt,0pt) .. controls (3pt,24pt)%
and (12pt,12pt) .. (24pt,0);

\draw[xshift=279pt,yshift=0pt] (0pt,0pt) .. controls (12pt,-12pt)%
and (21pt,-24pt) .. (24pt,0);
\node at(237pt,10pt){$\ldots$};
\draw[dotted] plot coordinates %
      {(0pt,0pt)(303pt,0pt)};
\draw plot coordinates %
      {(303pt,0pt)};

      \node at(70pt,40pt){\small{$M_{1}$}};
     \node at(148pt,20pt){\small{$M_{i-1}$}};
      \node at(200pt,10pt){\small{$M_{i}$}};
\node at(278pt,24pt){\small{$\overline{M_{k+1}}$}};
\node at (0pt,20pt){$P$};

\draw [xshift=350pt]plot coordinates %
  {(0,0) (10pt,0pt)(20pt,0pt)(30pt,0pt)(40pt,0pt)};
  \node at(370pt,-10pt){$H^{l}$};

\draw[very thick, color=blue,xshift=350pt,] plot coordinates %
  {(40pt,0pt)(50pt,-10pt)};
  \draw[xshift=400pt,yshift=-10pt] (0pt,0pt) .. controls (20pt,-20pt)%
and (35pt,-40pt) .. (40pt,0);
\draw[xshift=400pt,yshift=-10pt,very thick, color=blue] plot coordinates %
  {(40pt,0pt)(50pt,-10pt)};
 \draw[dashed,mark=,xshift=400pt,yshift=-10pt] plot coordinates %
      {(0pt,0pt)(40pt,0pt)};

      \draw[xshift=478pt,yshift=-30pt,very thick, color=blue] plot coordinates %
  {(-10pt,10pt)(0pt,0pt)};
       \draw[xshift=478pt,yshift=-30pt] (0pt,0pt) .. controls (20pt,-20pt)%
and (35pt,-40pt) .. (40pt,0);
  \draw[dashed,mark=,xshift=478pt,yshift=-30pt] plot coordinates %
      {(0pt,0pt)(40pt,0pt)};

 \draw[xshift=478pt,yshift=-30pt,very thick, color=red] plot coordinates %
  {(40pt,0pt)(50pt,50pt)};
\draw[xshift=528pt,yshift=20pt,very thick, color=blue] plot coordinates %
  {(40pt,0pt)(50pt,-10pt)};
  \draw[xshift=528pt,yshift=20pt] (0pt,0pt) .. controls (5pt,40pt)%
and (20pt,20pt) .. (40pt,0);
  \draw[dashed,mark=,xshift=528pt,yshift=20pt] plot coordinates %
      {(0pt,0pt)(40pt,0pt)};
\draw[dashed,xshift=350pt]plot coordinates %
      {(0pt,0pt)(40pt,0pt)};

      \node at(459pt,-20pt){$\ldots$};

 \draw[xshift=555pt,yshift=10pt,very thick, color=blue] plot coordinates %
  {(40pt,0pt)(50pt,-10pt)};
    \draw[xshift=605pt,yshift=0pt] (0pt,0pt) .. controls (3pt,24pt)%
and (12pt,12pt) .. (24pt,0);

\draw[xshift=629pt,yshift=0pt] (0pt,0pt) .. controls (12pt,-12pt)%
and (21pt,-24pt) .. (24pt,0);
\node at(587pt,10pt){$\ldots$};
\draw[dotted,xshift=350pt] plot coordinates %
      {(0pt,0pt)(303pt,0pt)};
\draw [xshift=350pt]plot coordinates %
      {(303pt,0pt)};

      \node at(422pt,1pt){\small{$\overline{M_{1}}$}};
     \node at(500pt,-18pt){\small{$\overline{M_{i-1}}$}};
      \node at(550pt,10pt){\small{$M_{i}$}};
\node at(625pt,24pt){\small{$\overline{M_{k+1}}$}};
\node at (350pt,20pt){$P_{i}$};
\end{tikzpicture}
\caption{An illustration of the decomposition of $P$ and the corresponding $P_i$.}\label{decomposeP}
\end{center}
\end{figure}
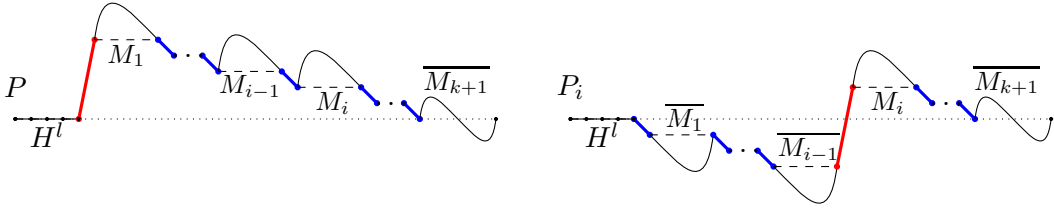

\begin{lem}\label{SPu} There is a 1-to-$(k+1)$ correspondence between $\mathcal {SP}^U_{n}(k,a)$ and $\mathcal{SP}^0_{n}(k,a)$, and we have
\begin{equation}\label{eqSPu}
|\mathcal{SP}_{n}^{U}(k,a)|=\frac{1}{k+1}\left(|\mathcal {SP}_{n}(k,a)|-\delta_{a|n}\right).
\end{equation}
\end{lem}
\pf
For each path $P$ in $\mathcal {SP}^U_{n}(k,a)$, we can uniquely decompose it into the following form:
$$P=H^l\ U\ M_1\ D\ M_2\ D\ \cdots\ D\ M_k\ D\ \overline{M_{k+1}},$$
in which $U$, $D$, and $H$ are single up, down, and horizontal steps,
respectively, $l\in\mathbb{N}$, $M_1,\ldots,M_k$ are $(k,a)$-paths, and
$\overline{M_{k+1}}$ is a super $(k,a)$-path (see Figure \ref{decomposeP}).


Now we map $P$ to the following $k+1$ paths in $\mathcal {SP}_{n}(k,a)$:
$$\psi(P)=\{P_i=H^l\ D\ \overline{M}_1\ D\ \cdots\ D\ \overline{M}_{i-1}\ {U}\ M_i\ D\ \cdots \ D\  \overline{M_{k+1}}: 1 \leqslant i \leqslant k+1\}.$$
Here $\overline{M}_i$ means the super $(k,a)$-path obtained from $M_i$ by reading the steps in reverse order, e.g., if $M_i=HUUDHD$, then $\overline{M}_i=DHDUUH$.  Therefore if $M_i$ is a $(k,a)$-path (never goes below the $x$-axis), then $\overline{M}_i$ is a super $(k,a)$-path that never goes above the $x$-axis. It is easy to see that $P_1=P$, and for $P_2, P_3,\ldots,P_{k+1}$, the first non-horizontal step is always a down step.

On the other hand, given any super $(k,a)$-path $P_{i}\in\mathcal{SP}^0_{n}(k,a)$, we can find a unique path $P\in\mathcal{SP}^U_{n}(k,a)$ such that $P_i\in \psi(P)$. First we find the left-most up step $U$ in $P_{i}$ whose right end point has non-negative $y$-coordinate, then decompose $P_i$ into the following form:
$$P_i=H^l\ D\ \overline{M_1}\ D\ \cdots\ D\ \overline{M_{i-1}}\ U\ M_i\ D\ \cdots \ D\  \overline{M_{k+1}}.$$
Then we have
$$P=H^l\ U\ M_1\ D\ \cdots\ D\ M_{i-1}\ D\ M_i\ D\ \cdots \ D\ \overline{M_{k+1}}.$$
Therefore, we proved that $\psi$ is a $1$-to-$(k+1)$ map from $\mathcal{SP}^U_{n}(k,a)$ to $\mathcal{SP}^0_{n}(k,a)$, and \eqref{eqSPu} follows. \qed

Now we are ready to give our main bijection in this section.

\begin{theo}\label{th:mainka} Let $\mathcal{LP}_{n}(k,a)$ denote the set of ordered pairs $(L,p)$, where $L\in \mathcal{P}_{n}(k,a)$, and $p$ is a specified hump in $L$. Then there is a bijection $\Phi:\mathcal{LP}_{n}(k,a)\rightarrow\mathcal{SP}_{n}^{U}(k,a)$.
\end{theo}
\pf Suppose $L\in \mathcal{P}_{n}(k,a)$, and $A$ is a lattice point on $L$. We use $x_{A}$ and $y_{A}$ to denote the $x$-coordinate and $y$-coordinate of $A$, respectively. The sub-path of $L$ from point $A$ to point $B$ is denoted by $L_{AB}$. If $p$ is a specified hump in $L$, we will also use the letter $p$ to denote the left end point of the first horizontal step in this specified hump and call it a \emph{hump point}. If the specified hump contains no horizontal steps, the hump point $p$ is the right end point of the up step in the hump.

For any $(L,p)\in\mathcal{LP}_{n}(k,a)$, we define the map $\Phi:\mathcal{LP}_{n}(k,a)\rightarrow\mathcal{SP}_{n}^{U}(k,a)$ as follows.
\begin{itemize}
\item  Let $A$ be the leftmost point in $L_{Op}$ that is followed by an up step, and there is no down step in $L_{Ap}$;
\item  Let $B$ be the leftmost point in $L$ such that $x_{B}>x_{p}$ and $y_{B}=y_{A}$;
\item  Let $C$ be the rightmost point in $L_{OA}$ such that $y_{C}=0$;
\item  Define $\Phi(L,p)=L_{OC} L_{AB}\overline{L}_{CA}\overline{L}_{BN}\triangleq SL$.
\end{itemize}

Since we only change the order of the steps, and the number of each type of step remains unchanged, we get a free $(k,a)$-path that ends at $(n,0)$. Moreover, the first non-horizontal step in $\Phi(L,p)$ is an up step: if $L_{OC}$ is empty or contains only horizontal steps, then the first non-horizontal step in $\Phi(L,p)$ is the up step following $A$, otherwise it is the first up step in $L_{OC}$. Therefore we proved that $\Phi(L,p)\in \mathcal {SP}_{n}^{U}(k,a)$.

Figure \ref{Fig1} shows, as an example, a $(3,1)$-path $L \in \mathcal{P}_{41}(3,1)$ with a specified hump point $p$, and Figure \ref{Fig2} shows the super $(3,1)$-path $\Phi(L,p) \in \mathcal{SP}_{41}^{U}$.

\begin{figure}[!htbp]
\centering
\begin{center}
\begin{picture}(400,110)
\put(-10,0){\circle*{4}} \put(0,0){\circle*{2}}
\put(10,30){\circle*{2}}\put(20,30){\circle*{2}}
\put(50,0){\circle*{2}}\put(60,0){\circle*{4}}
\put(70,30){\circle*{2}}\put(80,30){\circle*{2}}
\put(90,60){\circle*{2}} \put(100,50){\circle*{2}}\put(110,50){\circle*{4}}
\put(120,80){\circle*{2}}
\put(120,80){\circle{5}}
\put(130,80){\circle*{2}}
\put(140,70){\circle*{2}}
\put(150,100){\circle*{2}} \put(180,70){\circle*{2}}
\put(190,70){\circle*{2}} \put(210,50){\circle*{4}}
\put(220,50){\circle*{2}} \put(240,30){\circle*{2}}
\put(250,60){\circle*{2}} \put(260,60){\circle*{2}}
\put(310,10){\circle*{2}}
\put(320,10){\circle*{2}} \put(330,0){\circle*{2}}
\put(340,0){\circle*{2}} \put(350,30){\circle*{2}}
\put(380,0){\circle*{2}}
\put(400,0){\circle*{4}}

\put(30,20){\circle*{2}} \put(40,10){\circle*{2}}
\put(160,90){\circle*{2}} \put(170,80){\circle*{2}}
\put(200,60){\circle*{2}}
\put(230,40){\circle*{2}} \put(270,50){\circle*{2}}
\put(280,40){\circle*{2}} \put(290,30){\circle*{2}}
\put(300,20){\circle*{2}}
\put(360,20){\circle*{2}}\put(370,10){\circle*{2}}

\put(-10,0){\line(1,0){10}} \put(0,0){\line(1,3){10}}
\put(10,30){\line(1,0){10}}
\put(20,30){\line(1,-1){30}}
\put(50,0){\line(1,0){10}} \put(60,0){\line(1,3){10}}
\put(70,30){\line(1,0){10}} \put(80,30){\line(1,3){10}}
\put(90,60){\line(1,-1){10}}
\put(100,50){\line(1,0){10}}
\put(110,50){\line(1,3){10}}\put(120,80){\line(1,0){10}}
\put(130,80){\line(1,-1){10}}
\put(140,70){\line(1,3){10}}
\put(150,100){\line(1,-1){30}}
\put(180,70){\line(1,0){10}}
\put(190,70){\line(1,-1){20}}
\put(210,50){\line(1,0){10}}\put(220,50){\line(1,-1){20}}
\put(240,30){\line(1,3){10}} \put(250,60){\line(1,0){10}}
\put(260,60){\line(1,-1){50}} \put(310,10){\line(1,0){10}}
\put(320,10){\line(1,-1){10}} \put(330,0){\line(1,0){10}}
\put(340,0){\line(1,3){10}} \put(350,30){\line(1,-1){30}}
\put(380,0){\line(1,0){20}}

\put(117,84){\makebox(0,0)[bl]{p}}
\put(-14,-12){\makebox(0,0)[bl]{O}}
\put(106,38){\makebox(0,0)[bl]{A}}
\put(206,38){\makebox(0,0)[bl]{B}}
\put(56,-12){\makebox(0,0)[bl]{C}}
\put(396,-12){\makebox(0,0)[bl]{N}}
\multiput(-10,0)(0,10){12}{\qbezier[205](0,0)(205,0)(410,0)}
\multiput(-10,0)(10,0){42}{\qbezier[100](0,0)(0,50)(0,110) {\thinlines}}

\end{picture}
\end{center}
\caption{\label{Fig1}A $(3,1)$-path $L \in \mathcal{P}_{41}(3,1)$ with a circled hump point $p$.}
\end{figure}

\begin{figure}[!htbp]
\centering
\begin{center}
\begin{picture}(400,90)
\put(-10,10){\circle*{4}} \put(0,10){\circle*{2}}
\put(10,40){\circle*{2}}\put(20,40){\circle*{2}}
\put(50,10){\circle*{2}}\put(60,10){\circle*{4}}
\put(70,40){\circle{5}}
\put(70,40){\circle*{2}}\put(80,40){\circle*{2}}

 \put(90,30){\circle*{2}}
\put(100,60){\circle*{2}} \put(130,30){\circle*{2}}

\put(140,30){\circle*{2}}
\put(160,10){\circle*{4}}
\put(170,10){\circle*{2}} \put(180,0){\circle*{2}}
\put(190,30){\circle*{2}} \put(200,30){\circle*{2}}
\put(210,60){\circle*{4}} \put(230,60){\circle*{2}}
\put(260,30){\circle*{2}} \put(270,60){\circle*{2}}
\put(280,60){\circle*{2}}
\put(290,50){\circle*{2}} \put(300,50){\circle*{2}}
\put(350,0){\circle*{2}} \put(360,0){\circle*{2}}
\put(370,30){\circle*{2}} \put(390,10){\circle*{2}}
\put(400,10){\circle*{4}}

\put(30,30){\circle*{2}} \put(40,20){\circle*{2}}
\put(110,50){\circle*{2}} \put(120,40){\circle*{2}}
\put(240,50){\circle*{2}}
\put(250,40){\circle*{2}} \put(310,40){\circle*{2}}
\put(320,30){\circle*{2}} \put(330,20){\circle*{2}}
\put(340,10){\circle*{2}} \put(380,20){\circle*{2}}

\put(-10,10){\line(1,0){10}} \put(0,10){\line(1,3){10}}
\put(10,40){\line(1,0){10}}
\put(20,40){\line(1,-1){30}}
\put(50,10){\line(1,0){10}} \put(60,10){\line(1,3){10}}
\put(70,40){\line(1,0){10}}

\put(80,40){\line(1,-1){10}} \put(90,30){\line(1,3){10}}
\put(100,60){\line(1,-1){30}}
\put(130,30){\line(1,0){10}}
\put(140,30){\line(1,-1){20}}\put(160,10){\line(1,0){10}}
\put(170,10){\line(1,-1){10}}
\put(180,0){\line(1,3){10}}
\put(190,30){\line(1,0){10}}
\put(200,30){\line(1,3){10}}
\put(210,60){\line(1,0){20}}

\put(230,60){\line(1,-1){30}}\put(260,30){\line(1,3){10}}
\put(270,60){\line(1,0){10}} \put(280,60){\line(1,-1){10}}
\put(290,50){\line(1,0){10}} \put(300,50){\line(1,-1){50}}
\put(350,0){\line(1,0){10}} \put(360,0){\line(1,3){10}}
\put(370,30){\line(1,-1){20}} \put(390,10){\line(1,0){10}}
\multiput(-10,0)(0,10){9}{\qbezier[205](0,0)(205,0)(410,0)}
\multiput(-10,0)(10,0){42}{\qbezier[80](0,0)(0,50)(0,80) }

\qbezier[205](0,10)(205,10)(400,10)

\put(60,40){\makebox(0,0)[bl]{p}}
\put(-14,-4){\makebox(0,0)[bl]{O}}
\put(56,-4){\makebox(0,0)[bl]{A}}
\put(156,-4){\makebox(0,0)[bl]{B}}
\put(200,60){\makebox(0,0)[bl]{C}}
\put(396,-4){\makebox(0,0)[bl]{N}}
\end{picture}
\end{center}
\caption{\label{Fig2}A super $(3,1)$-path $\Phi(L,p) \in \mathcal{SP}^{U}_{41}(3,1)$.}
\end{figure}

Now we define the inverse of $\Phi$. For any super $(k,a)$-path $SL\in \mathcal{SP}_{n}^{U}(k,a)$:

\begin{itemize}
\item Let $B$ be the point on the $x$-axis that follows a down step, and the next down step in $SL$ is the first down step that goes below the $x$-axis. If $SL$ is always above the $x$-axis, then we set $B$ to be the end point $N$.
\item  Let $A$ be the rightmost point with $y_A=0$, $x_A < x_B$ and $A$ is followed by an up step.
\item  Let $C$ be the leftmost point in $SL$ such that $ x_{C}\geqslant x_{B}$ and $\forall G$, $x_{G}\geqslant x_{B}$ implies that $y_{G}\leqslant y_{C}$;
\item  Let $p$ be the leftmost hump in $L_{AB}$, and set $L=L_{OA}\overline{L}_{BC}L_{AB}\overline{L}_{CN}$;
\item Set $\Psi(SL)=(L,p)$.
\end{itemize}
It is easy to check that $\Psi=\Phi^{-1}$. Therefore $\Phi$ is a bijection. \qed

\noindent{\it Proof of Equation \eqref{hump} and \eqref{peak}}:
\eqref{hump} follows immediately from Lemma \ref{SPu} and Theorem \ref{th:mainka}.

For Equation \eqref{peak}, we will prove it by counting the number of humps in all $(k,a)$-paths in $\mathcal{P}_n(k,a)$ that are not peaks. Given $(L,p)\in\mathcal{LP}_{n}(k,a)$, if the specified hump $p$ in $L$ is not a peak, then in the resulting super $(k,a)$-path $SL=\Phi(L,p)=L_{OC} L_{AB}\overline{L}_{CA}\overline{L}_{BN}$, the leftmost hump in $L_{AB}$ is not a peak. If we remove the first horizontal step in this hump in $SL$, we get a super $(k,a)$-path $\widetilde{SL} \in \mathcal{SP}_{n-a}^U(k,a)$. On the other hand, given any super $(k,a)$-path $\widetilde{SL}\in \mathcal{SP}_{n-a}^U(k,a)$, we first find the essential points $A$ and $B$ in $\widetilde{SL}$ as we define the map $\Psi$ in the proof of Theorem \ref{th:mainka}, and then add a horizontal step to the first hump in $L_{AB}$ of $\widetilde{SL}$ and get $SL \in \mathcal{SP}_{n}^U(k,a)$. Set $(L,p)=\Psi(SL)$. Then $p$ is a hump in $L$ that is not a peak. Hence we've established a bijection between $\mathcal{SP}^U_{n-a}(k,a)$ and the set of ordered pairs $(L,p)\in\mathcal{LP}_{n}(k,a)$ such that $p$ is not a peak.

Therefore by applying Lemma \ref{SPu} we have that the total number of peaks in all $(k,a)$-paths of order $n$ is
\begin{eqnarray*}
\sum_{P\in\mathcal{P}_{n}(k,a)}\#\Peaks(P)&=&\frac{1}{k+1}\left(|\mathcal {SP}_{n}(k,a)|-\delta_{a|n}\right)-
\frac{1}{k+1}\left(|\mathcal {SP}_{n-a}(k,a)|-\delta_{a|(n-a)}\right)\\
&=&\frac{1}{k+1}(|\mathcal{SP}_{n}(k,a)|-|\mathcal{SP}_{n-a}(k,a)|).
\end{eqnarray*}
\qed

\noindent{\emph Remark 1:}
Yan also gave bijective proofs of \eqref{hump} and \eqref{peak} in \cite{SYan}  but her bijection is different from our bijection $\Phi$.

\noindent{\emph Remark 2:}
Note that when defining the bijection $\Phi$, the parameters $k$ and $a$ do not really matter. Let $S$ be a set of positive integers. We define an $(S,a)$-path of order $n$ to be a lattice path in $\mathbb{Z}\times\mathbb{Z}$ from $(0,0)$ to $(n,0)$ which uses up steps $U=(1,k), k\in S$, down steps $D=(1,-1)$ and horizontal steps $H=(a,0)$ and never goes below the $x$-axis. Therefore, our bijection $\Phi$ proves the following stronger result for $(S,a)$-paths:

\begin{coro}
The total number of humps in all $(S,a)$-paths of order $n$ equals the total number of super $(S,a)$-paths of order $n$ whose first non-horizontal step is an up step.
\end{coro}

\section{$(n,m)$-Dyck Paths with a given number of peaks}

Let $\mathcal{D}(n,m)$ and $\mathcal{F}(n,m)$ denote the set of $(n,m)$-Dyck paths and the set of free paths using $(1,0)$ and $(0,1)$ steps from $(0,0)$ to $(n,m)$, respectively. For any free path $P \in \mathcal{F}(n,m)$, we define the equivalence class of $P$ to be the set of all cyclic permutations of the steps making up $P$ that lead to distinct free paths from $(0,0)$ to $(n,m)$, and denote it as $[P]$. More precisely, suppose $P=u_{1}u_{2}\cdots u_{n+m}$, where $u_i\in \{U,D\}$ for each $i, 1\leq i \leq n+m$, then the equivalence class of $P$ is
\[[P]=\{P_i:= u_{i+1}u_{i+2}\cdots u_{n+m}u_1u_2\cdots u_i | i = 1, 2, \cdots,n+m\}.\]
Here $P_{n+m}=P$. For example, if $P=DUUDU$, then
\[[P]=\{P_1=UUDUD, P_2=UDUDU, P_3=DUDUU, P_4=UDUUD, P_5=P=DUUDU\}.\]

When $\gcd(n,m)=1$, the equivalence class of $P$ has the following properties.

\begin{lem}\label{lem-fukukawa}
For any free path $P$ from $(0,0)$ to $(n,m)$, if $\gcd(n,m)=1$, then
\begin{itemize}
\item[1)] $|[P]|=n+m$;
\item[2)] There is a unique $(n,m)$-Dyck path in $[P]$.
\end{itemize}
\end{lem}

\noindent
\pf \begin{itemize}
\item[1)]
For any path $P =u_{1}u_{2}\cdots u_{n+m} \in \mathcal{F}(n,m)$, let $r$ be the smallest positive integer such that $P_{r}=u_{r+1}u_{r+2}\cdots u_{n+m}u_{1}u_{2}\cdots u_{r}=P$. Then we have $|[P]|=r$. It is obvious that $r\leqslant n+m$. Now we claim that $r=n+m$. Otherwise, we have $P=P_{r}=P_{2r}=P_{3r}=\cdots=P_{n+m}$, therefore $n+m=ir$ for some positive integer $i\geq 2$. Suppose there are $x$ $D$'s and $y$ $U$'s in $u_{1}u_{2}\cdots u_{r}$. Then $n=ix$ and $m=iy$, which contradicts the condition that $\gcd(n,m)=1.$ Therefore we have $|[P]|=n+m$.

\item[2)]
For any path $P\in \mathcal{F}(n,m)$, if $P$ is not an $(n,m)$-Dyck path, then there must be at least one lattice point on $P$ that is below the diagonal line $y=\frac{m}{n}x$. Let $v$ be one of these points that is furthest away from the diagonal line. Then we can decompose $P$ into two sub-paths at the point $v$, i.e., $P=L_1 L_2$. Setting $\tilde{P}=L_2 L_1$, it is obvious that $\tilde{P} \in [P]$ (see Figure \ref{PtildeP}). Suppose $v$ is on the line $y=\frac{m}{n}x-c$ for some positive number $c$, then all the other lattice points on $P$ are above this line. Therefore in $\tilde{P}$ both $L_2$ and $L_1$ are above the diagonal line. Hence $\tilde{P}$ is an $(n,m)$-Dyck path .

\begin{figure}[h]
\begin{center}
\begin{tikzpicture}[scale=0.6,very thick,domain=0:8,mark=]
   \draw[thin,gray,step=15pt] (0,0) grid (75pt,120pt);
   \draw[yshift=0cm,blue] plot coordinates %
   {(0pt,0pt) (0pt,30pt) (30pt,30pt) (30pt,45pt) (45pt,45pt)};
   \draw[yshift=0cm] plot coordinates %
   {(45pt,45pt) (45pt,90pt) (60pt,90pt) (60pt,120pt) (75pt,120pt)};
  \draw[dashed,thin,gray] (0,0)--(75pt,120pt);
  \fill[red] (45pt,45pt) circle (2pt);
   \node at (53pt,45pt) [scale=0.8] {$v$};
   \node at (20pt,85pt) [scale=0.8] {$P$};
   \node at (25pt,20pt) [scale=0.8] {$L_{1}$};
   \node at (60pt,80pt) [scale=0.8] {$L_{2}$};
  \draw[xshift=5cm,thin,gray,step=15pt] (0,0) grid (75pt,120pt);
   \draw[xshift=5cm,yshift=0cm] plot coordinates %
   {(0pt,0pt) (0pt,45pt) (15pt,45pt) (15pt,75pt) (30pt,75pt) };
   \draw[xshift=5cm,yshift=0cm,blue] plot coordinates %
   {(30pt,75pt) (30pt,105pt) (60pt,105pt) (60pt,120pt) (75pt,120pt)};
   \fill[red,xshift=5cm] (30pt,75pt) circle (2pt);
   \draw[xshift=5cm,dashed,thin,gray] (0,0)--(75pt,120pt);
   \node at (20pt,95pt) [scale=0.8,xshift=4.2cm] {$\bar{P}$};
   \node at (10pt,30pt) [scale=0.8,xshift=4.2cm] {$L_{2}$};
   \node at (60pt,100pt) [scale=0.8,xshift=4.2cm] {$L_{1}$};
  \draw[xshift=2.3cm,yshift=1.2cm,->,thick] (15pt,0)--(70pt,0pt);
\end{tikzpicture}
\end{center}\caption{$P$ and $\tilde{P}$}\label{PtildeP}
\end{figure}
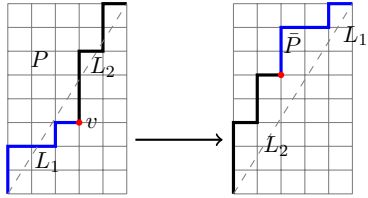

Now we will show that such a $\tilde{P}$ is unique. We claim that there are no two different points that are both furthest away from the diagonal line. Suppose there are two different points $v_1=(x_1, y_1)$ and $v_2=(x_2,y_2)$. Then the line connecting $v_1$ and $v_2$ is parallel to the line $y=\frac{m}{n}x$. Then we have
$$
\frac{y_{2}-y_{1}}{x_{2}-x_{1}}=\frac{m}{n},
$$
which contradicts the condition $\gcd(n,m)=1.$

Finally we want to show that for any $Q \in [P]$, $\tilde{P}=\tilde{Q}$. Suppose $\tilde{P} \neq \tilde{Q}$. Since $\tilde{P}$ and $\tilde{Q}$ are both $(n,m)$-Dyck paths, and they are in the same equivalence class $[P]$, there must be at least one lattice point $(x_0,y_0)$ on $\tilde{P}$ ($\tilde{Q}$) that lies on the diagonal line $y=\frac{m}{n}x$ with $0< x_0 < n$. But this is impossible when $\gcd(n,m)=1$.

Therefore we proved that there is a unique $(n,m)$-Dyck path in $[P]$ when $\gcd(n,m)=1$.
\end{itemize}
\qed

\begin{figure}[h]
\begin{center}
\begin{tikzpicture}[thick,domain=0:8,mark=]
\draw[thin,gray,step=15pt] (0,0) grid (30pt,45pt);
\draw[yshift=0cm,very thick] plot coordinates %
{(0pt,0pt) (15pt,0pt) (15pt,30pt) (30pt,30pt) (30pt,45pt)};
\coordinate [label=0:$P$] ($P$) at (8pt,-10pt);
\draw[thin,gray,xshift=4 cm,step=15pt] (0,0) grid (30pt,45pt);
\draw[xshift=4 cm,yshift=0cm,very thick] plot coordinates %
{(0pt,0pt) (0pt,30pt) (15pt,30pt) (15pt,45pt) (30pt,45pt)};
\coordinate [label=0:$P_{1}$] ($P_{1}$) at (122pt,-10pt);
\draw[thin,gray,xshift=6 cm,step=15pt] (0,0) grid (30pt,45pt);
\draw[xshift=6 cm,yshift=0cm,very thick] plot coordinates %
{(0pt,0pt) (0pt,15pt) (15pt,15pt) (15pt,30pt) (30pt,30pt) (30pt,45pt)};
\coordinate [label=0:$P_{2}$] ($P_{2}$) at (179pt,-10pt);
\draw[thin,gray,xshift=8 cm,step=15pt] (0,0) grid (30pt,45pt);
\draw[xshift=8 cm,yshift=0cm,very thick] plot coordinates %
{(0pt,0pt) (15pt,0pt) (15pt,15pt) (30pt,15pt) (30pt,45pt)};
\coordinate [label=0:$P_{3}$] ($P_{3}$) at (236pt,-10pt);
\draw[thin,gray,xshift=10 cm,step=15pt] (0,0) grid (30pt,45pt);
\draw[xshift=10 cm,yshift=0cm,very thick] plot coordinates %
{(0pt,0pt) (0pt,15pt) (15pt,15pt) (15pt,45pt) (30pt,45pt)};
\coordinate [label=0:$P_{4}$] ($P_{4}$) at (293pt,-10pt);
\draw[thin,gray,xshift=12 cm,step=15pt] (0,0) grid (30pt,45pt);
\draw[xshift=12 cm,yshift=0cm,very thick] plot coordinates %
{(0pt,0pt) (15pt,0pt) (15pt,30pt) (30pt,30pt) (30pt,45pt)};
\coordinate [label=0:$P_{5}$] ($P_{5}$) at (350pt,-10pt);
\draw[xshift=1.5 cm,yshift=0.5cm,->,thick] (0,0)--(50pt,0pt);
\end{tikzpicture}
\end{center}
\caption{A free path $P$ from $(0,0)$ to $(2,3)$, and the 5 different free paths in $[P]$.}\label{Pclass}
\end{figure}
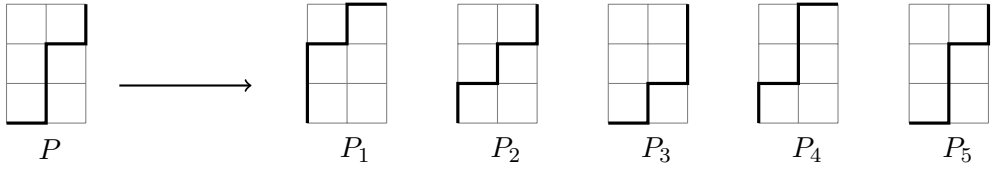

\begin{exam}\label{equ}
Figure \ref{Pclass} shows a free path $P$ from $(0,0)$ to $(2,3)$, and the $5$ different free paths in $[P]$, in which $P_1$ is the unique $(2,3)$-Dyck path.
\end{exam}

Since the total number of free paths from $(0,0)$ to $(n,m)$ is ${m+n \choose n}$, we immediately deduce the following result from Lemma \ref{lem-fukukawa}.
\begin{coro}\label{gc}
When $\gcd(n,m)=1$, then the number of $(n,m)$-Dyck paths is
\begin{equation}
D(n,m)=\frac{1}{n+m}{n+m\choose n}.
\end{equation}
\end{coro}

We remark that Lemma \ref{lem-fukukawa} and Corollary \ref{gc} are also proved in \cite{Bizley,Y.FUKUKAWA} in even stronger forms. Here we are more interested in the refined enumeration of $(n,m)$-Dyck paths with a given number of peaks and the relation with free $(n,m)$-paths.

Let $\mathcal{D}(n,m;j)$ ($\mathcal{F}(n,m;j)$) denote the set of $(n,m)$-Dyck paths (free paths) with exactly $j$ peaks. And let $\mathcal{F}^{UD}(n,m;j)$ denote the set of paths in $\mathcal{F}(n,m;j)$ that start with an up step and end with a down step. The following theorem shows that there are similar relations between peaks in all $(n,m)$-Dyck paths and $\mathcal{F}^{UD}(n,m;j)$ to what we found for $(k,a)$-paths in the previous section.

\begin{theo}\label{mainbij-nm}
Let $\mathcal{PD}(n,m;j)=\{(P,p)|P\in\mathcal{D}(n,m;j),\ p\mbox{ is a peak of } P\}$. Then there is a bijection $\hat{\Phi}:\mathcal{PD}(n,m;j)\rightarrow\mathcal{F}^{UD}(n,m;j)$ when $\gcd(n,m)=1$.
\end{theo}
\noindent\pf
Given $(P,p)\in\mathcal{PD}(n,m;j)$, since there are $j$ peaks in $P$, we can uniquely decompose $P$ into $P=L_{1}L_{2}\cdots L_{j}$, such that each $L_{i}$ is of the form $U^{a_i} D^{b_i}$ for some positive integers $a_i$ and $b_i$ since an $(n,m)$-Dyck path must start with a $U$ and end with a $D$.  Suppose $p$ is contained in $L_{i} (1\leqslant i \leqslant j)$. Then we set
$$
\hat{P}=\hat{\Phi}(P,p)=L_{i+1}L_{i+2}\cdots L_{j}L_{1}L_{2}\cdots L_{i}.
$$
Note that when $i=j$, we set  $\hat{P}=\hat{\Phi}(P,p)=P$.

It is easy to check that $\hat{P}\in\mathcal{F}^{UD}(n,m;j)$. Now we will show that $\hat{\Phi}$ is a bijection by defining the inverse of $\hat{\Phi}$. For any path $\hat{P}\in\mathcal{F}^{UD}(n,m;j)$, we can uniquely decompose it as $\hat{P}=L_{1}L_{2}\cdots L_{j}$, in which each $L_{i}$ is of the form $U^{a_i} D^{b_i}$ for some positive integers $a_i$ and $b_i$. Let $v$ be the point that is furthest below the line $y=\frac{m}{n}x$ among all the lattice points on $\hat{P}$. Then $v$ must be the point that connects $L_i$ and $L_{i+1}$ for some $i, 1\leqslant i\leqslant j-1$. Moreover, from the proof of Lemma \ref{lem-fukukawa} we know that such a $v$ is unique. Set $P=L_{i+1}L_{i+2}\cdots L_{j}L_{1}L_{2}\cdots L_{i}$, and let $p$ be the only peak in $L_j$. Since $\gcd(n,m)=1$, from the proof of Lemma \ref{lem-fukukawa} we know that $P$ is the unique $(n,m)$-Dyck path in $[\hat{P}]$. Let $(P,p)=\hat{\Psi}(\hat{P})$, then we have $(P,p)\in\mathcal{PD}(n,m;j)$ and $\hat{\Psi}=\hat{\Phi}^{-1}$. Therefore we proved that $\hat{\Phi}$ is a bijection. \qed

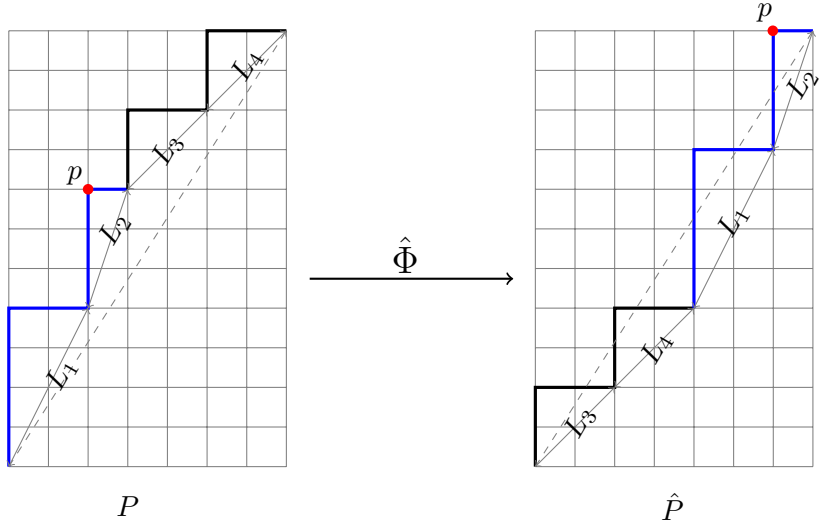
\begin{figure}[h]
\begin{center}
\begin{tikzpicture}[very thick,domain=0:8,mark=]
   \draw[thin,gray,step=15pt] (0,0) grid (105pt,165pt);
   \draw[yshift=0cm,blue] plot coordinates %
   {(0pt,0pt) (0pt,60pt) (30pt,60pt) (30pt,105pt) (45pt,105pt)};
   \draw[yshift=0cm] plot coordinates %
   {(45pt,105pt) (45pt,135pt) (75pt,135pt) (75pt,165pt) (105pt,165pt)};
  \fill[red] (30pt,105pt) circle (2pt);
  \draw[dashed,thin,gray] (0,0)--(105pt,165pt);
  \draw[thin,gray,<->] (0,0)--(30pt,60pt);
  \draw[thin,gray,<->] (30pt,60pt)--(45pt,105pt);
  \draw[thin,gray,<->] (45pt,105pt)--(75pt,135pt);
  \draw[thin,gray,<->] (75pt,135pt)--(105pt,165pt);
  \node at (25pt,110pt)  {$p$};
  \node at (45pt,-15pt)  {$P$};
  \node at (20pt,35pt) [rotate=60] {$L_{1}$};
  \node at (40pt,90pt) [rotate=60] {$L_{2}$};
  \node at (60pt,120pt) [rotate=60] {$L_{3}$};
  \node at (90pt,152pt) [rotate=60] {$L_{4}$};
  \draw[xshift=7cm,thin,gray,step=15pt] (0,0) grid (105pt,165pt);
   \draw[xshift=7cm,yshift=0cm] plot coordinates %
   {(0pt,0pt) (0pt,30pt) (30pt,30pt) (30pt,60pt) (60pt,60pt)};
   \draw[xshift=7cm,yshift=0cm,blue] plot coordinates %
   {(60pt,60pt) (60pt,120pt) (90pt,120pt) (90pt,165pt) (105pt,165pt)};
   \fill[red] (289pt,165pt) circle (2pt);
   \draw[xshift=7cm,dashed,thin,gray] (0,0)--(105pt,165pt);
  \draw[xshift=7cm,thin,gray,<->] (0,0)--(30pt,30pt);
  \draw[xshift=7cm,thin,gray,<->] (30pt,30pt)--(60pt,60pt);
  \draw[xshift=7cm,thin,gray,<->] (60pt,60pt)--(90pt,120pt);
  \draw[xshift=7cm,thin,gray,<->] (90pt,120pt)--(105pt,165pt);
  \node at (286pt,172pt)  {$p$};
  \node at (251pt,-15pt)  {$\hat{P}$};
  \node at (216pt,17pt) [rotate=60] {$L_{3}$};
  \node at (245pt,45pt) [rotate=60] {$L_{4}$};
  \node at (274pt,93pt) [rotate=60] {$L_{1}$};
  \node at (300pt,145pt) [rotate=60] {$L_{2}$};
  \draw[xshift=4cm,yshift=2.5cm,->,thick] (0,0)--(77pt,0pt);
  \node at (150pt,80pt) [rotate=0] {\Large{$\hat{\Phi}$}};
\end{tikzpicture}
\end{center}
\caption{A $(7,11)$-Dyck path $P$ with a specified peak $p$, and the coresponding free path $\hat{P}$.}\label{gduck}
\end{figure}

\begin{exam}
As an example, Figure \ref{gduck} shows a $(7,11)$-Dyck path $P$ with a specified peak $p$, and the coresponding free path $\hat{P}=\hat{\Phi}(P,p)$.
\end{exam}

From Theorem \ref{mainbij-nm} we know that if we can find the number of free paths in $\mathcal{F}^{UD}(n,m;j)$, we will be able to find the number of $(n,m)$-Dyck paths with exactly $j$ peaks. And the following Lemma counts the number of free paths in $\mathcal{F}^{UD}(n,m;j)$.

\begin{lem}
\begin{itemize}
\item[1)] The number of free paths from $(0,0)$ to $(n,m)$ with $j$ peaks is
\begin{equation}\label{gpeak}
|\mathcal{F}(n,m;j)|={n\choose j}{m\choose j};
\end{equation}
\item[2)] The number of free paths from $(0,0)$ to $(n,m)$ with $j$ peaks that start with an up step and end with a down step, is
\begin{equation}\label{UDgpeak}
|\mathcal{F}^{UD}(n,m;j)|={n-1\choose j-1}{m-1\choose j-1}.
\end{equation}
\end{itemize}
\end{lem}

\noindent\pf
1)Let $P \in \mathcal{F}(n,m;j)$. For each peak $UD$ in $P$, we call the lattice point between $U$ and $D$ a \emph{peak point}. It is obvious that a free path with $j$ peaks is uniquely determined by the $j$ peak points $p_i=(x_i,y_i), i=1,\ldots, j$, with $0\leqslant x_{1}<x_{2}<\cdots<x_{j}\leqslant n-1$, and $1\leqslant y_{1}<y_{2}<\cdots<y_{j}\leqslant m$, and this is a one-to-one correspondence. We have ${n\choose j}$ ways to choose the $x_i$'s, and ${m \choose j}$ ways to choose the $y_i$'s. Therefore there are ${n\choose j}{m\choose j}$ ways to choose these peak points. Hence we proved Equation \eqref{gpeak}.

2)If $P\in \mathcal{F}^{UD}(n,m;j)$, the coordinates $(x_{1},y_{1}),(x_{2},y_{2}),\cdots, (x_{j},y_{j})$ of the $j$ peak points of $P$ must satisfy $0=x_{1}<x_{2}<\cdots<x_{j}\leqslant n-1$ and $1\leqslant y_{1}<y_{2}<\cdots<y_{j}=m$. Therefore we can choose the $j$ peak points in ${n-1\choose j-1}{m-1\choose j-1}$ ways. Hence we get \eqref{UDgpeak}. \qed

Combining Theorem \ref{mainbij-nm} and Equation \eqref{UDgpeak}, we get the following theorem.
\begin{theo}\label{dnm}
When $\gcd(n,m)=1$, the number of $(n,m)$-Dyck paths with exactly $j$ peaks is:
\begin{equation}\label{rationalNRYN}
D(n,m;j)=\frac{1}{j}{n-1\choose j-1}{m-1\choose j-1}.
\end{equation}
\end{theo}

We remark here that \eqref{rationalNRYN} is also given in \cite{DArmstrong}, in which the authors call it a ``rational Narayana number".

\section{$k$-ary paths with a given number of peaks}

By taking special values of $m$ for the results we obtained in the previous section, we can get the number of $k$-ary paths with a given number of peaks.

\begin{lem}\label{1-1-1}
There are one-to-one correspondences between the following sets:  $(n,kn)$-Dyck paths, $(n,kn+1)$-Dyck paths, and $k$-ary paths of order $(k+1)n$.
\end{lem}
\pf For each path $P\in \mathcal{D}(n,kn+1)$, the first two steps of $P$ are both up steps. (Otherwise, it will go below the diagonal line). By deleting the first up step, we get a path $P^{'}$ that goes from $(0,1)$ to $(n,kn+1)$, and never goes below the line that connects $(0,1)$ and $(n,kn+1)$ (no lattice point lies above the line $y=\frac{kn+1}{n}x$ and below the line $y=kx+1$ for $x \leq n$.). It is obvious that these paths are in one-to-one correspondence with $(n,kn)$-Dyck paths. Moreover, for each $(n,kn)$-Dyck path $P^{'}$, if we take the steps of $P^{'}$ in reverse order and then swap the $U$'s and $D$'s, we get a unique $k$-ary path $P^{''}$ of order $(k+1)n$. (Here we regard $P^{'}$ and $P^{''}$ as sequences of $U$'s and $D$'s). It is easy to see that such a correspondence is also one-to-one.\qed

\begin{figure}[h]
\begin{center}

\begin{tikzpicture}[scale=0.6,domain=0:8,mark size=0.12ex]
   \draw[thin,gray,step=15pt] (0,0) grid (75pt,165pt);
   \draw[yshift=0cm,very thick,red] plot coordinates %
   {(0pt,0pt) (0pt,15pt)};
   \draw[yshift=0cm,very thick] plot coordinates %
  {(0pt,15pt) (0pt,60pt) (15pt,60pt) (15pt,135pt) (45pt,135pt) (45pt,165pt) (75pt,165pt)};
  \draw[thin,gray] (-5pt,-11pt)--(80pt,176pt);
  \fill[red] (0pt,0pt)  circle (2pt);
  \fill[red] (0pt,15pt) circle (2pt);
  \fill[red] (75pt,165pt) circle (2pt);
     \node at (80pt,120pt) [scale=0.8] {$y=\frac{kn+1}{n}x$};
   \node at (40pt,-15pt)  {$P$};
    \node at (-18pt,15pt) [scale=0.8] {$(0,1)$};
     \node at (-15pt,-8pt) [scale=0.8] {$(0,0)$};
      \node at (85pt,175pt) [scale=0.8] {$(n,kn+1)$};
   \draw[thin,xshift=6cm,gray,step=15pt] (0,0) grid (75pt,165pt);
   \draw[xshift=6cm,very thick] plot coordinates %
   {(0pt,0pt) (0pt,45pt) (15pt,45pt) (15pt,120pt) (45pt,120pt) (45pt,150pt) (75pt,150pt)};
   \draw[thin,gray,xshift=6cm] (-5pt,-10pt)--(80pt,160pt);
    \node at (80pt,100pt)[scale=0.8,xshift=3.6cm]  {$y=kx$};
   \node at (40pt,-15pt)[xshift=3.6cm]  {$P'$};
    \fill[red] (0pt,0pt) [xshift=6cm] circle (2pt);
     \fill[red] (75pt,150pt)[xshift=6cm] circle (2pt);
      \node at (18pt,0pt) [scale=0.8,xshift=3.6cm] {$(0,0)$};
      \node at (113pt,160pt) [scale=0.8,xshift=3.6cm] {$(n,kn)$};

 \draw[thin,gray,step=15pt,xshift=320pt] (0pt,0pt) grid (225pt,105pt);
   \draw[yshift=0cm,thick, xshift=320pt,mark=*] plot coordinates %
   {(0pt,0pt) (15pt,30pt)(30pt,60pt) (45pt,45pt) (60pt,30pt) (75pt,60pt)
   (90pt,90pt) (105pt,75pt) (120pt,60pt)(135pt,45pt) (150pt,30pt) (165pt,15pt)
   (180pt,45pt) (195pt,30pt)(210pt,15pt) (225pt,0pt)};
   \fill[red] (320pt,0pt)  circle (2pt);
   \fill[red] (545pt,0pt)  circle (2pt);
    \node at (0pt,-15pt)[xshift=260pt]  {$P''$};
\end{tikzpicture}
\end{center}
\caption{\label{fig7} A $(5,11)$-Dyck path, a $(5,10)$-Dyck path and a $2$-ary path of order $15$.}
\end{figure}
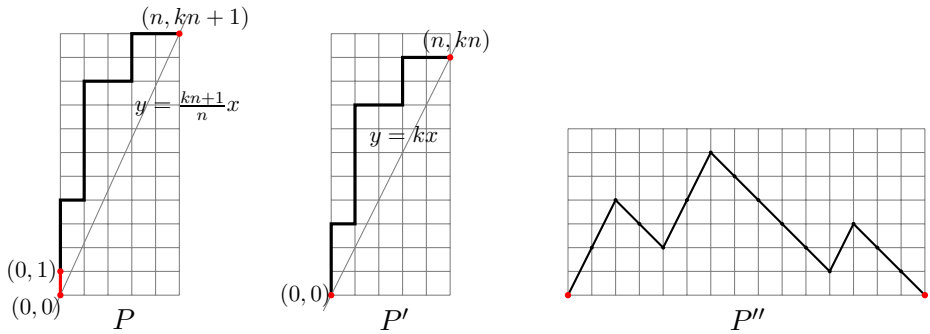

Figure \ref{fig7} shows an example of a $(5,11)$-Dyck path, the corresponding $(5,10)$-Dyck path and the corresponding $2$-ary path of order $15$.

For any integer $n$, if $m=kn+1$ for some positive integer $k$, we always have that $\gcd(n,m)=\gcd(n,kn+1)=1$.
Therefore from Lemma \ref{1-1-1}, Corollary \ref{gc} and Theorem \ref{dnm} we have the following results on $k$-ary paths.

\begin{coro}
\begin{itemize}
\item The number of $k$-ary paths of order $(k+1)n$ is:
\begin{equation}
\frac{1}{kn+1}{(k+1)n\choose n}; \label{dnkn}
\end{equation}
\item The number of $k$-ary paths of order $(k+1)n$ with exactly $j$ peaks is:
\begin{equation}D(n,k;j)=\frac{1}{j}{n-1\choose j-1}{kn\choose j-1}. \label{dnknj}
\end{equation}
\end{itemize}
\end{coro}

\noindent {\em Remark:} Note that when $k=1$, equation \eqref{dnkn} and \eqref{dnknj} coincide with the well-known result that Dyck paths of order $n$ are counted by the $n$-th Catalan number
$C(n)=\frac{1}{n+1}{2n \choose n}$, and the number of Dyck paths of order $n$ with exactly $j$ peaks is the Narayana number $N(n;j)=\frac{1}{j}{n-1\choose j-1}{n\choose j-1}$. And counting peaks of height $k$ in a Dyck path has been studied by Mansour in \cite{Mansour_peakheightk}.

\vskip 2mm \noindent{\bf Acknowledgments.} The authors would like to thank the referees for carefully reading an earlier version of this paper and for giving many constructive comments which helped improving the quality of the paper a lot. This work is partially supported by the National Science Foundation of China under Grant No. 10801053, Shanghai Rising-Star Program (No. 10QA1401900), and the Fundamental Research Funds for the Central Universities.

\end{document}